\documentclass[12pt]{amsart}
\usepackage[left=1in,top=1in,right=1in,bottom=1in,letterpaper]{geometry}
\usepackage{epsfig}
\usepackage{amsmath}
\usepackage{amssymb}
\usepackage{amscd}
\usepackage{bm}
\usepackage{url}
\usepackage{verbatim} 
\usepackage{color}
\usepackage{epsfig}
\usepackage{stmaryrd}
\usepackage{amsthm}
\usepackage{pifont}
\usepackage{mathrsfs}
\usepackage{wasysym}

\usepackage{hyperref}

\title{Generalized Staircases: Recurrence and Symmetry}
\author[W. Patrick Hooper]{W. Patrick Hooper$^*$}
\thanks{$^*$ Supported by N.S.F. Postdoctoral Fellowship DMS-0803013}
\address{Department of Mathematics, Northwestern University\\
2033 Sheridan Road\\
Evanston, IL 60208-2730, USA (phone: 847-491-2853, fax: 847-491-8906)}
\email{wphooper@math.northwestern.edu}
\author[Barak Weiss]{Barak Weiss$^{\dagger}$}
\thanks{$^{\dagger}$ Supported by BSF grant 2004149 and ISF grant 584/04}
\address{Ben Gurion University, Be'er Sheva, Israel 84105}
\email{ barakw@math.bgu.ac.il}

\newtheorem{theorem}{Theorem}
\newtheorem{proposition}[theorem]{Proposition}
\newtheorem{lemma}[theorem]{Lemma}
\newtheorem{remark}[theorem]{Remark}
\newtheorem{corollary}[theorem]{Corollary}

\theoremstyle{definition}
\newtheorem{definition}[theorem]{Definition}

\def\P{\mathbb{P}}%
\def\Q{\mathbb{Q}}%
\def\R{\mathbb{R}}%
\def\Z{\mathbb{Z}}%
\def\GL{\textit{GL}}
\def\SL{\textit{SL}}

\def\PGL{\textit{PGL}}

\renewcommand{\hom}[1]{\ensuremath{{\llbracket #1 \rrbracket}}}%
\def\ker{\textit{ker}}%
\newcommand{\nullset}{\emptyset}

\def\Aut{\textit{Aut}}%
\def\H{\mathbb H}%

\def\isomS1{\textit{Isom}_+(S^1)}
\def\rt3{\sqrt{3}}

\def\tr{\mathrm{tr}} 

\def\Isom{\mathrm{Isom}}

\newif\ifdraft\drafttrue
\draftfalse

\long\def\combarak#1{\ifdraft{\tt #1 }\else\ignorespaces\fi}
\long\def\compat#1{\ifdraft{\textcolor{blue}{\tt #1 }}\else\ignorespaces\fi}

\def\H{{\mathbb H}} 

\def\Aff{{\textit A\hspace{-1pt}f\hspace{-3pt}f}}

\def\rk{{\mathrm{rk}}}
\def\dim{{\mathit{dim} }}
\def\hol{{\mathbf{hol}}}                                 
\def\Fix{{\mathit{Fix}}}

\begin{document}
\maketitle
\combarak{I don't know why, but the figures and the references
misbehave when I compile the file on my computer.}
\compat{I compile this with "pdflatex Zcovers" (maybe twice) then add references with "bibtex Zcovers" then run "pdflatex Zcovers" again. I hope this helps...}
\begin{abstract}
We study infinite translation surfaces which are $\Z$-covers of
compact translation surfaces. We obtain conditions ensuring that such
surfaces have Veech groups which are Fuchsian of the first
kind and give a necessary and sufficient condition for recurrence of
their straight-line flows. Extending results of Hubert and Schmith\"usen, we
provide examples of infinite non-arithmetic lattice surfaces, as well
as surfaces with infinitely generated Veech groups. 
\end{abstract}
\section{Introduction}
The geometry of translation surfaces has been intensively studied in
recent years. While most of the work was concerned with compact
surfaces, in several recent papers non-compact surfaces were also
considered. For instance, in \cite{CGL}, the horseshoe and baker's transformations
were realized by an affine transformation; \cite{Higl} is a study of the geometry and dynamics of an infinite translation
surfaces which arises as a geometric limits of compact lattice surfaces; in \cite{HW08}, a connection was made to $\Z$-valued skew
products over 1-dimensional systems, and in \cite{Valdez}, the
topology of the unfolding surface for an irrational billiard was
determined. Removing the restriction that the surface is compact gives
a flexible setup and many phenomena, absent in the compact case, may
be observed. For example, in the recent paper \cite{HS09}, Hubert and
Schmith\"usen made the surprising discovery that there are infinite 
square tiled surfaces whose Veech group is infinitely generated. 

The examples studied in \cite{HW08, HS09} are $\Z$-covers of compact
translation surfaces. Although this class is much smaller than the
general case, it already displays many surprising features. It may be
hoped that it provides a good starting point 
for a study of the geometry and dynamics of infinite translation 
surfaces. In this paper we begin the systematic study of these
surfaces. Our analysis yields a bijection between $\Z$-covers
$\widetilde M \to M$, ramified over a finite set $P \subset M$,  and
projective classes of elements $w \in H_1(M, P; \Z)$ (Proposition
\ref{prop 5}). Under this bijection, {\em recurrent
$\Z$-covers}, i.e. covers on which the straightline flow is recurrent
in almost all directions, correspond to homology classes with
vanishing holonomy (Proposition \ref{prop: recurrence}). Utilizing a
result of Thurston \cite{Thurston98}, we obtain a
sufficient condition ensuring that Veech groups of a cover $\widetilde
M$ is Fuchsian of the first kind (Theorem \ref{thm:dimension_2}). This result implies
that any $\Z$-cover of a square tiled surface in genus 2 has a Veech group which is of the
first kind (Corollary \ref{cor2}), extending the results of
\cite{HS08}. \compat{Note change in wording below. (old version "We also obtain necessary and sufficient conditions for
a (finite power of a) parabolic element in the Veech group of $M$ to lift to the Veech
group of $\widetilde M$ (Theorem \ref{thm:multi-twist})."} We also obtain necessary and sufficient conditions for
a (finite power of a) parabolic element in the Veech group of $M$ to lift to the Veech
group of every recurrent $\Z$-cover $\widetilde M$ (Theorem \ref{thm:multi-twist}). Using it one
may reprove some of the results of \cite{HS08} in a more general
setting. We illustrate the use of our 
results in the last section, where we provide an example of an
infinite lattice translation surface with a non-arithmetic Veech
group (Proposition \ref{prop: octagon}), and answer a question of Hubert and
Schmith\"usen (Section \ref{sect:HS}).
\section{Regular covers of translation surfaces}

Let $M$ denote a compact translation surface and $P \subset M$ denote
a finite (possibly empty) subset. We consider $P$ to be a collection
of punctures of the surface $M$ and will use $M^\circ$ to denote $M
\smallsetminus P$. 

Covering space theory associates covers of a space with the subgroups of its fundamental group. 
A cover is called {\em regular} if it is associated to a normal subgroup.
We consider a normal subgroup $N \subset \pi_1(M^\circ)$,
and consider the associated cover $\pi:\widetilde M \to M^\circ$. 
The group $\Delta = \pi_1(M^\circ)/N$ acts on $\widetilde M$ as the
automorphisms of the cover, with $\widetilde M/\Delta=M^\circ$.  

We have the following from covering space theory.
\begin{proposition}~\vspace{-0.25em}
\label{prop:covering_theory}
\begin{enumerate}
\item An element $f \in \Aff(M^\circ)$ lifts to an $\widetilde f \in
\Aff(\widetilde M)$ if and only if $f_\ast(N)=N$. 
\item An element $\widetilde f \in \Aff(\widetilde M)$ descends to an $f \in
\Aff(M^\circ)$ if and only if $\widetilde f$ normalizes the deck group
$\Delta$. 
That is, $\widetilde f \Delta \widetilde f^{-1}=\Delta$. 
\end{enumerate}
\end{proposition}

\begin{definition} The {\em affine automorphism group of a cover
$\widetilde M \to M^\circ$} is the group of pairs of elements 
$( \widetilde f, f) \in \Aff(\widetilde M) \times \Aff(M^\circ)$ for which
$\pi \circ \widetilde f=f \circ \pi$. We denote this group by $\Aff(\widetilde
M, M^\circ)$. A necessary condition for 
$( \widetilde f, f) \in \Aff(\widetilde M, M^\circ)$ is that $D(\widetilde
f)=D(f)$. Thus we have a canonical definition of the derivative
$D:\Aff(\widetilde M, M^\circ) \to \GL(2, \R)$. We call the image 
of the group homomorphism $D$ the {\em Veech group of the cover}, and
denote it by $\Gamma(\widetilde M)$. 
\end{definition}

Let $G_N=\{f \in \Aff(M^\circ)~:~f_\ast(N)=N\}$. For an $f \in G_N$
the action of $f_\ast$ on $\pi_1(M^\circ)$ induces 
an action on $\Delta=\pi_1(M^\circ)/N$. The following is an immediate consequence:

\begin{corollary}
$\Aff(\widetilde M, M^\circ) \cong \Delta \rtimes G_N$, with $G_N$ acting on $\Delta$ as mentioned above.
Indeed, we have a short exact sequence
$$1 \to \Delta \hookrightarrow \Aff(\widetilde M, M^\circ) \twoheadrightarrow G_N \to 1.$$
\end{corollary}
\qed

Note that the projection $p:\Aff(\widetilde M, M) \to \Aff(\widetilde M)$ may
not be injective. However, we do not miss much. 

\begin{proposition}
If $M^\circ$ is not an unpunctured torus,
then $p\big(\Aff(\widetilde M, M^\circ)\big)$ is a finite index subgroup of $\Aff(\widetilde M)$. 
\end{proposition}
\begin{proof}
Consider the group ${\mathcal T} \subset \Aff(\widetilde M)$ of elements
$\iota$ for which $D(\iota)=I$, i.e. the group of translation
automorphisms of $\widetilde M$. We claim that ${\mathcal T}$ acts properly discontinuously on the set of
non-singular points of $\widetilde M$. To see this, let $Q$ denote the union of the singularities of $M$ with $P$. By assumption $Q$ is non-empty.
The surface $M^\circ$ has a Delaunay decomposition relative to the points in $Q$. See \cite[\S 1]{MS91} for background.
The Delaunay decomposition of $\widetilde M$ relative to the lifts of $Q$ is the lift of this decomposition to $M^\circ$. 
A translation automorphism must permute the cells in the Delaunay decomposition, and hence is properly discontinuous.

The deck group $\Delta$ is a finite index subgroup of ${\mathcal T}$, because
$\mathrm{Area}(\widetilde M/\Delta)=\mathrm{Area}(M)<\infty$.
The group $\Delta$ is finitely generated because it is a quotient of $\pi_1(M^\circ)$, which is finitely generated.
${\mathcal T}$ is also finitely generated as it contains $\Delta$ as a finite index subgroup. An element
$\widetilde f \in \Aff(\widetilde M)$ acts on ${\mathcal T}$ by 
conjugation, and preserves the index of subgroups. There are only
finitely many subgroups of ${\mathcal T}$ with index 
$[{\mathcal T}:\Delta]$, because ${\mathcal T}$ is finitely generated.
Thus, a finite index subgroup of $\Aff(\widetilde M)$ normalizes $\Delta$. The conclusion follows by 
Proposition \ref{prop:covering_theory}.
\end{proof}

\combarak{This begs the question of what happens when $\Delta$ is not
finitely generated. Do you know?}
\compat{I realized that $\Delta$ is always finitely generated! See the modified proof above.}
\section{$\Z$-covers}

We use $H_1(M, P; \Z)$ to denote the relative homology of $M$ with
respect to the set of punctures, and $H_1(M^\circ; \Z)$ denotes the
absolute homology of 
the punctured surface. Intersection number is a non-degenerate bilinear form 
$$i:H_1(M, P; \Z) \times H_1(M^\circ; \Z) \to \Z.$$

\begin{definition}
The {\em $\Z$-cover of $M^\circ$ associated to a non-zero $w \in
H_1(M, P; \Z)$} is the cover associated to the kernel of the
homomorphism 
$$\varphi_{w}:\pi_1(M^\circ) \to \Z, \ \ \gamma \mapsto i(w,
\hom{\gamma}),$$
where $\hom{\gamma}$ denotes the homology class of $\gamma$. 
We denote this cover by $\widetilde M_w$. 
\end{definition}

If $A$ is a free abelian group, we use $\P A$ to denote 
$(A \smallsetminus \{{\mathbf 0}\})/\sim,$
where $a \sim b$ if there are non-zero $m,n \in \Z$  for which $ma=nb$. 
The discussion above shows:
\begin{proposition}
\label{prop 5}
The $\Z$-covers $\widetilde M_{w}$ and $\widetilde M_{w'}$ are the same
if and only if $w \sim w'$. 
\end{proposition}

Thus, the space of $\Z$-covers of $M^\circ$ is naturally identified with $\P H_1(M, P; \Z)$.
Statement (1) of proposition \ref{prop:covering_theory} can be restated as follows.

\begin{proposition}
\label{prop:aff}
An $f \in \Aff(M^\circ)$ lifts to an $\widetilde f \in \Aff(\widetilde M_w)$ if and only if $f_\ast(w)=\pm w$. 
\end{proposition}

A translation surface has a {\em holonomy map} $\hol:H_1(M, P; \Z) \to
\R^2$, obtained by developing a representative of the class into $\R^2$ and taking the difference of the starting and end points. 
For dynamical reasons, we are especially interested in
$\Z$-covers with the following property. 

\begin{definition}[Recurrent $\Z$-covers]
The $\Z$-cover, $\widetilde M_w$ is called {\em recurrent} if $\hol(w)={\mathbf 0}$. 
\end{definition}

Although not explicitly stated, square-tiled covers of this type were
studied before in \cite{HW08} and \cite{HS09}. One reason for
restricting attention to  
recurrent $\Z$-covers is that non-recurrent $\Z$-covers have few
affine symmetries. A subgroup of $\GL(2, \R)$ is called 
{\em elementary} if it contains a finite index abelian subgroup and
{\em non-elementary} otherwise. Conversely, a non-elementary subgroup
of $\GL(2, \R)$ 
contains the free group with two generators. We have the following corollary of Proposition \ref{prop:aff}.

\begin{corollary}
If $\Gamma(\widetilde M_w)$ is non-elementary, then $\widetilde M_w$ is a recurrent $\Z$-cover.
\end{corollary}
\begin{proof}
We prove the contrapositive. Suppose $\hol(w) \neq {\mathbf 0}$. By
Proposition \ref{prop:aff}, $(\widetilde f, f) \in (\widetilde M_w, M^\circ)$
implies 
that $f_\ast(w)=\pm w$. Then
$D(f)\big(\hol(w)\big)=\hol\big(f_\ast(w)\big)=\pm \hol(w).$
In particular any element of $\Gamma(\widetilde M_w)$ has $\hol(w)$ as an eigenvector.
\end{proof}

We will now justify the term recurrent $\Z$-cover. Let $F_t^{\theta}:M
\to M$ denote the straight-line flow in direction $\theta \in S^1$.  
Similarly, we will use $\widetilde F_t^{\theta}:\widetilde M \to \widetilde M$ to
denote the straight-line flow on a $\Z$-cover $\widetilde M$ in direction
$\theta$. Recall that a flow is measure preserving flow $F_t$ is
called {\em recurrent} if for any measurable set $A$, for a.e. $x \in
A$ there is $t_n \to \infty$ such that $F_{t_n}x \in A$. 

\begin{proposition}[Recurrence of the straight-line flow]\label{prop: recurrence}
Let $\widetilde M$ be a $\Z$-cover of $M^\circ$. 
Then $\widetilde M$ is a recurrent $\Z$-cover if and only if for any $\theta$ for which
$F^{\theta}_t$ is ergodic, $\widetilde F^{\theta}_t$ is recurrent.
\end{proposition}

\begin{proof}
We will reduce the statement to a classical result of K. Schmidt
\cite[Theorem 11.4]{Schmidt} in infinite ergodic theory. 
Suppose $(X,
\mu)$ is a finite measure 
space and $T: X \to X$ is a measurable transformation preserving $\mu$
which is ergodic. 
For a measurable $f: X \to \Z$, $f \in L^1(X, \mu)$, define $X_f = X
\times \Z$ and  
$$T_f : X_f \to X_f, \ \ T_f(x,k) = \left(Tx, k+f(x) \right).$$
Then $T_f$ is recurrent if
and only if $\int f \, d\mu=0$. 

\compat{Minor modifications of proof. We want $\widetilde \alpha$ to be the set of all lifts of $\alpha$, rather than just one lift. Also note that
$\alpha$ cannot be disjoint from $w$, as $w$ is just a homology class, but I do not think this affects the proof.}

\compat{{\bf Old version}:
Given $\theta$, we reduce to the above statement as follows: 
choose segments $\alpha, \widetilde
\alpha$ in $M, \widetilde M$ respectively, which are in the direction
$\theta'$ perpendicular to $\theta$, such that $\widetilde
\alpha$ is a lift of $\alpha$, and such that $\alpha$ is disjoint from
$w$.}

Given $\theta$, we reduce to the above statement as follows: 
choose a segment $\alpha$ in $M$, which is the  direction
$\theta'$ perpendicular to $\theta$. Define $\widetilde \alpha$ in
$\widetilde M$ to be the union of all lifts of $\alpha$ to $\widetilde M$. 
Denote by $T$ (resp. $\widetilde T$) the Poincar\'e return map to
the section $\alpha$ (resp. $\widetilde \alpha$), so that $T$ is an
interval exchange transformation. The ergodicity of
$T$ is equivalent to that of $F^{\theta}_t$ and the recurrence of $\widetilde
F_t^{\theta}$ is equivalent to the 
recurrence of $\widetilde T$. Since continuous
maps have Borel sections, we may (measurably) identify $\widetilde M$
with $M \times \Z$. In these coordinates $\widetilde T  = T_f$ where
$$f=f^{(\theta)}: \alpha \to \Z, \ \ f(x) = i(w, \hom{\beta_x}),$$
and $\beta_x$ is the 
curve from $x$ to $Tx$ along the $F^{\theta}_t$ orbit of $x$, and then
from $Tx$ to $x$ along $\alpha$. Let $\mu$ be the length measure
on $\alpha$. Up to scaling, Lebesgue measure on $M$ can be 
represented as $d\mu dt$, where $dt$ denotes the length measure along the orbits of
$F^{\theta}_t$. Since $f$ assumes finitely
many values, one on each interval of continuity of $T$, it is in
$L^1(\alpha, 
\mu)$.

Label by $I_1, 
\ldots, I_\ell$ be the partition of $\alpha$ into intervals of
continuity for $T$. By refining this 
decomposition we assume that the flow in direction $\theta$ starting
from the interior of $I_j$ does not hit a puncture in $P$. For each $j$, let
$\beta_j$ be a closed loop $\beta_{x_j}$ as above, corresponding to
some $x_j \in I_j$; the particular choice of $x_j$ does not affect
$\hom{\beta_j}$. 
Now write $\beta = \sum \mu(I_j)
\hom{\beta_j}  \in H_1(M; \R)$. We claim that for a path $\gamma$ on $M$
representing an element of $H_1(M, P; \Z)$, 
$$i(\gamma, \hom{\beta}) = \hol_{\theta'} (\gamma),$$
i.e. the holonomy in the direction $\theta'$. Indeed after homotoping $\gamma$ off $\alpha$, each
positive crossing of $\ell_j$ means $\gamma$ has crossed the rectangle
above $I_j$, and contributes $\mu(I_j)$ to
$\hol_{\theta'}(\gamma)$. Therefore 
\[
\begin{split}
\int f^{(\theta)} \, d\mu & = \sum_j \mu(I_j) i(w, \hom{\gamma_j}) \\
& = i(w, \hom{\beta}) 
 = \hol_{\theta'}(w). 
\end{split}
\]
The main theorem of \cite{KMS85} guarantees the existence of two independent ergodic $\theta$. We
see that
$\int f^{(\theta)} \, d\mu =0$ for any ergodic direction $\theta$ on $M$, if and
only if $\hol (w)=0$. 
\end{proof}
\combarak{I vaguely remember that in Tel
Aviv your argument for the multitwist theorem involved a calculation
similar to the one above, relating $\hol_{\theta'}(w)$ to an
intersection number, and now I don't see any such calculation. Maybe the
proof can be shortened, or maybe the calculation can contribute to the
proof of the multitwist theorem and make it shorter, what do you think?}
\compat{I think that was my original argument, but I found it hard to make it completely rigorous. I think the current proof is better.}

\begin{corollary}
If $\hol(w)=0$, the straightline flow $\widetilde F^{\theta}_t$ on $\widetilde M_w$ is recurrent for a.e. $\theta$.
\end{corollary}
\begin{proof}
Combine Proposition \ref{prop: recurrence} with the famous result of Kerckhoff, Masur and Smillie \cite{KMS85}.
\end{proof}

\section{Veech groups of recurrent $\Z$-covers}

Let $H \subset \R^2$. We define $K(H)$ to be the smallest extension
field of $\Q$ for which there is an $A \in \GL(2, \R)$ such that $A(H)
\subset K(H)^2$.  
As in the appendix of \cite{KS}, the {\em holonomy field} of a
translation surface $M$ is the field $k=K\big(\hol(H_1(M; \Z))\big)$.  
Kenyon and Smillie showed that if $M$ is compact and there is a
pseudo-Anosov in $\Aff(M)$, then $k$ is a field 
extension of $\Q$ of degree at most the genus $g$ of $M$. Moreover,
the image $\hol\big(H_1(M; \Z)\big)$ is a $\Z$-module of rank
$2[k:\Q]$.  

It follows from the work of Kenyon and Smillie that if $\Aff(M^\circ)$
contains a pseudo-Anosov then  
$K\big(\hol(H_1(M;\Z))\big)=K\big(\hol(H_1(M,P;\Z))\big)$. We
unambiguously declare this the holonomy field in this case, and 
we use $k$ to denote this field.

\begin{definition}[Holonomy-free subspaces]
The {\em holonomy-free subspaces of homology} are $W=\ker~\hol \subset
H_1(M, P; \Z)$ of relative homology, 
and $W_0=W \cap H_1(M;\Z)$ of absolute homology.
\end{definition}

The $\Z$-modules $W_0$ and $W$ have ranks given by the following equations.
$$\rk~W_0=\rk~H_1(M; \Z)-2[k:\Q]=2(g-[k:\Q]).$$
$$\rk~W=\rk~H_1(M,P; \Z)-2[k:\Q]=\begin{cases}
2(g-[k:\Q])+\#P-1 & \text{if $P \neq \nullset$}\\
2(g-[k:\Q]) & \textrm{otherwise.}
\end{cases}
$$ 

The affine automorphism group $\Aff(M^\circ)$ acts on homology and
preserves the subspaces $W_0$ and $W$. Thus, we have the following
group homomorphisms. 
$$\psi_0:\Aff(M^\circ) \to \Aut(W_0), \ \ f \mapsto f_\ast|_{W_{0}}.$$
$$\psi:\Aff(M^\circ) \to \Aut(W), \ \ f \mapsto f_\ast|_{W}.$$
The following statement follows immediately from Proposition 
\ref{prop:aff}. It 
explains our interest in these homomorphisms.

\begin{proposition}
Let $f \in \ker~\psi$. For each $w \in W$, there is an $\widetilde f \in
\Aff(\widetilde M_w)$ such that $(\widetilde f, f) \in \Aff(\widetilde M_w,
M^\circ)$. The subgroup 
$$\{ (\widetilde f, f) \in \Aff(\widetilde M_w, M^\circ) ~:~ f \in \ker~\psi \}$$
is normal inside $\Aff(\widetilde M_w, M^\circ)$.
\end{proposition}

The elements of $\Aff(M^\circ)$ permute the punctures. 
Let $\rho:\Aff(M^\circ) \to \mathit{Sym}(P)$ be the map which assigns
to an $f \in \Aff(M^\circ)$ the permutation induced 
on $P$. We have the following. 

\begin{proposition}
\label{prop:kernels}
$\psi(\ker~\psi_0 \cap \ker~\rho)$ is abelian of rank at most
$(\rk~W_0)(\rk~W-\rk~W_0)$. Thus, there is an exact sequence 
$$1\to \ker~\psi \hookrightarrow \ker~\psi_0 \twoheadrightarrow A \to 1$$
where $A \subset \Z^{(\rk~W_0)(\rk~W-\rk~W_0)} \rtimes
\mathit{Sym}(P)$ has a finite index free abelian subgroup. 
\end{proposition}
\begin{proof}
Enumerate $P=\{p_1, \ldots, p_n\}$, and let $\gamma_i \in H_1(M^\circ;
\Z)$ be the homology class of a loop which travels clockwise around
$p_i$ 
for $i=1,\ldots, n$. Let $J:W \to \Z^n$ denote the function
$$J(w)=\big( i(w,\gamma_1), \ldots,  i(w,\gamma_n)\big)\in \Z^n.$$ 
Note that for all $f  \in \Aff(M^\circ)$ we have $J \circ
f_\ast(w)=\rho(f) \circ J(w)$, where the permutation $\rho(f)$ is
acting as a permutation matrix. 
In addition, $J(w)$ determines the coset of $W/W_0$ which contains
$w$. The following statements follow from this discussion. 
\begin{enumerate}
\item $\ker~J=W_0$.
\item If $f \in \ker~\rho$, then $f_\ast(w)-w \in W_0$ for all $w \in W$.  
\end{enumerate}
By definition, if $f \in \ker~\psi_0$, then $f_\ast(w_0)=w_0$ for all $w_0 \in W_0$. 
For $f \in \ker~\psi_0 \cap \ker~\rho$, let $h_f:W/W_0 \to W_0$ denote
the map $w+W_0 \mapsto f_\ast(w)-w$. 
This is well defined by the above discussion.
Moreover, we can recover $\psi(f)=f_\ast|_W$ via the formula
$\psi(f)(w)=w+h_f(w+W_0)$. If $f, g \in \ker~\psi_0 \cap \ker~\rho$,
$$\begin{array}{rcl}
\psi(g \circ f)(w)=\psi(g)\big(w+h_{f}(w+W_0)\big) & = & w+h_{f}(w+W_0)+h_g\big(w+h_{f}(w+W_0)+W_0\Big)\\
& = & w+h_f(w+W_0)+h_g(w+W_0).
\end{array}$$
So $\psi(\ker~\psi_0 \cap \ker~\rho)$ is abelian group. Moreover, an
element $\psi(f)$ of this group is uniquely determined by  
the linear map $h_f:W/W_0 \to W_0$. It can be observed that 
$W/W_0 \cong \Z^{\rk~W-\rk~W_0}$ and $W_0 \cong \Z^{\rk~W_0}$. Hence, the space of
all possible $h_f$ is isomorphic to $\Z^{(\rk~W_0)(\rk~W-\rk~W_0)}$.
\end{proof}

If $G$ is a discrete subgroup of $\GL(2, \R)$, we will use $\Lambda G
\subset \R \P^1$ to denote the limit set of the projection of $G$ 
to $\PGL(2, \R)=\Isom(\H^2)$. 
A subgroup $G$ of $\GL(2, \R)$ or $\PGL(2, \R)$ is elementary if and
only if $\Lambda G$ contains two or fewer points.  
See \cite{MT98} for background on the limit set and for the following.

\begin{lemma}[Limit sets of normal subgroups]
\label{lem:normal}
Suppose $G$ is a non-elementary discrete subgroup of $\GL(2, \R)$ or
$\PGL(2, \R)$. If $N$ is a non-trivial normal subgroup of $G$, then
$\Lambda N=\Lambda G$. 
\end{lemma}

\begin{theorem}
\label{thm:kernel}
If $D\big(\Aff(M^\circ)\big)$ is non-elementary and $D(\ker~\psi_0)$ is non-trivial, then 
$$\Lambda D\big(\Aff(M^\circ)\big)=\Lambda D(\ker~\psi_0)=\Lambda D(\ker~\psi).$$
In this case, $\Lambda \Gamma (\widetilde M_w)=\Lambda D\big(\Aff(M^\circ)\big)$
for all recurrent $\Z$-covers $\widetilde M_w$ of $M^\circ$. 
\end{theorem}
\begin{proof}
If $D(\ker~\psi_0)$ is non-trivial, then by a direct application of Lemma \ref{lem:normal},
$\Lambda D\big(\Aff(M^\circ)\big)=\Lambda D(\ker~\psi_0)$. In particular, $D(\ker~\psi_0)$
contains a free group with two generators.
By Proposition \ref{prop:kernels},
$D(\ker~\psi)$ is a finite index subgroup of the kernel of a map from $D(\ker~\psi_0)$
to an abelian group. Hence, $D(\ker~\psi)$ is non-empty. By another application of Lemma \ref{lem:normal},
we see $\Lambda D(\ker~\psi)=\Lambda D(\ker~\psi_0)$.
\end{proof}

A {\em Fuchsian group of the first kind} is a subgroup $\Gamma$ of
$\Isom(\H^2)$ (or some other linear group 
which acts isometrically on $\H^2$) for which $\Lambda \Gamma=\R \P^1$. 

\begin{theorem}
\label{thm:dimension_2}
Suppose $D\big(\Aff(M^\circ)\big)$ is a lattice and that $\rk~W_0 \leq 2$.
Then $D(\ker~\psi)$ is a Fuchsian group of the first kind. In particular, for any $w \in W$,
$\Gamma (\widetilde M_w)$ is Fuchsian of the first kind.
\end{theorem}
\begin{proof}
By Theorem \ref{thm:kernel}, it is sufficient to show that $D(\ker~\psi_0)$ is non-trivial. 
Note that $\rk~W_0$ is even.
If $\rk~W_0=0$, then $\ker~\psi_0=\Aff(M^\circ)$. The more difficult
case is when $\rk~W_0=2$. We will assume that 
$\ker~\psi_0$ is empty and derive a contradiction.

As observed by Veech \cite{V}, $D\big(\Aff(M^\circ)\big)$ is not co-compact. Thus
$D\big(\Aff(M^\circ)\big)$ contains a free group of finite rank, which
pulls back to a free group $F \subset \Aff(M^\circ)$ such that $D|_F$ is injective. 
Since $\rk~W_0=2$, $\psi_0:F \to \hat \SL(2, \Z)$, where $\hat \SL(2,
\Z)$ denotes the set of $2 \times 2$ matrices 
of determinant $\pm1$. By assumption, $\psi_0|_F$ is injective. 
Without loss of generality, we may assume that $\psi_0(F) \subset
\SL(2, \Z)$. (If not, replace $F$ by the index two subgroup for which
this is true.) 

Summarizing the previous paragraph, we have two faithful
representations, $D|_F$ and $\psi_0|_F$, of $F$ into $\SL(2, \R)$. We
will derive a contradiction from properties 
of these representations. These representations satisfy the following statements for all $f \in F$. 
\begin{enumerate}
\item If $D(f)$ is parabolic, then $\psi_0(f)$ is also parabolic.
\item If $D(f)$ is hyperbolic, then $2 \leq |\tr~\psi_0(f)| < |\tr~D(f)|$. 
\end{enumerate}
Statement 1 is true because if $f \in \Aff(M^\circ)$ is a parabolic,
then some power of $f$ is a multi-twist of $M^\circ$. All eigenvalues
of the action of a multi-twist 
on homology are $1$. In particular, the eigenvalues for the action of
$f$ on homology are all of modulus $1$. Thus, $\psi_0(f)$ is either
elliptic or parabolic. But,  
if $\psi_0(f)$ is elliptic, then $\psi_0$ is not faithful. If $D(f)$
is hyperbolic, then $f \in \Aff(M^\circ)$ is a pseudo-Anosov. Let
$\lambda$ be the eigenvalue 
of $D(f)$ with largest magnitude. A theorem of Fried implies that
$\lambda$ is also the eigenvalue with largest magnitude 
of the action of $f_\ast$ on $H_1(M^\circ; \Z)$, and also that
$\lambda$ occurs with multiplicity one \cite{Fried85}. 
In particular, the eigenvalues of $\psi_0(f)=f_\ast|_{W_0}$ have modulus strictly less than $|\lambda|$. 
Again, $\psi_0(f)$ is not elliptic since $\psi_0$ is assumed to be faithful.

Now consider the quotient surfaces
$S_1=\H^2/D(F)$ and $S_2=\H^2/\psi_0(F)$.
For $i=1,2$, let $g_i$ denote the genus of $S_i$ and let $n_i \geq 1$ denote the number of ends. 
We have $F=\pi_1(S_1)=\pi_1(S_2)$, so this induces a homotopy equivalence $\phi:S_1 \to S_2$. 
Thus, we have that $\rk~F=2 g_i+n_i-1$ for each $i$. By statement $1$ above, we have $n_1 \leq n_2$.
We will show that $g_1=g_2$ and $n_1=n_2$. 

An element of the fundamental group of a surface is called {\em
peripheral} if it is homotopic to a puncture. 
Assume that $n_1<n_2$. 
Let $\gamma_1, \ldots, \gamma_{n_1} \in \pi_1(S_1)$ denote disjoint
peripheral curves. Note that the homology classes 
of these curves are linearly dependent. Let $\gamma_j'=\phi_\ast(\gamma_j) \in \pi_1(S_2)$. Note that since 
$S_2$ has $n_2>n_1$ punctures, the homology classes of the curves
$\gamma_1', \ldots, \gamma_{n_1}'$ are linearly independent. 
This contradicts either the fact that $\phi$ is a homotopy
equivalence, or that $n_1<n_2$. Thus, $n_1=n_2$.  

By the previous two paragraphs, we may take the homotopy equivalence
$\phi:S_1 \to S_2$ to be a homeomorphism. In addition, these surfaces
have the same number of parabolic cusps. 
Thus $\psi_0(F)$ is a lattice in $\SL(2, \Z)$. For non-peripheral
$\beta \in \pi_1(S_1)$ let $\ell_1(\beta)$ denote the length of the
geodesic representative 
on $S_1$, and let $\ell_2(\beta)$ denote the length of the geodesic
representative of $\phi_\ast(\beta)$. Theorem 3.1 of \cite{Thurston98}
states that 
$$\sup_{\beta \in \pi_1(S_1)} \frac{\ell_2(\beta)}{\ell_1(\beta)} \geq 1,$$
with equality only if $S_1=S_2$. (This holds for any pair of complete,
finite area, hyperbolic structures on the same surface.) This
contradicts statement (2).
\end{proof}

The following immediate consequence illustrates the use of Theorem
\ref{thm:dimension_2}. 

\begin{corollary}
\label{cor2}
If $M$ is any translation surface of genus 1 or 2 with non-elementary
Veech group, then the Veech group of any recurrent $\Z$-cover has the
same limit set. In particular, if $M$ is a square tiled surface of
genus 1 or 2 then the Veech group of any recurrent $\Z$-cover is
Fuchsian of the first kind. 

\end{corollary}

\section{Multi-twists}

A {\em multi-twist} is an $f \in \Aff(M^\circ)$
which preserves the cylinders in a cylinder decomposition and for
which $D(f)$ is parabolic with eigenvalue $1$. 
It is well known that if $M$ is compact, and $D(f)$ is parabolic then
some power of $f$ is a multi-twist. The action of a multi-twist $f$ on
$H_1(M,P;\Z)$ is given by the formula 
\begin{equation}
\label{eq:action_on_homology}
f_\ast:x \mapsto x+\sum_{j} i(x, \gamma^\circ_j) t_j \gamma_j,
\end{equation}
where $j$ varies over the cylinders in the preserved decomposition. Here $\gamma_j \in H_1(M, P; \Z)$ 
and $\gamma^\circ_j \in H_1(M^\circ;\Z)$ denote
the homology classes of the core curve in cylinder
$j$ (although the curves are the same they represent
elements in different homology spaces, and we will use different
notation to distinguish them). We denote by $\langle \gamma_j
\rangle$ and $\langle \gamma^\circ_j \rangle$ the $\Z$-module spanned
by these curves in their respective homology groups. The
restriction of the action of $f$ on cylinder $j$ is a Dehn twist. The
number $t_j \in \Z$ is the twist number of this Dehn twist. 
Each $t_j$ is non-zero and they all have the same sign. If this sign
is positive $f$ is performing left Dehn twists and if it is negative
$f$ is performing right Dehn twists. 

Let $\phi=f_\ast-I$. That is, 
\begin{equation}
\label{eq:phi}
\phi:H_1(M,P;\Z) \to H_1(M,P;\Z), \ \ x \mapsto \sum_{j} i(x, \gamma^\circ_j) t_j \gamma_j.
\end{equation}
A direct application of Proposition \ref{prop:aff} yields the following. \compat{While nearly trivial, I feel this is important enough to state as a proposition:}

\begin{proposition}
The multi-twist $f \in \Aff(M^\circ)$ lifts to an $\widetilde f \in \Aff(\widetilde M_w, M^\circ)$ if and only if
$\phi(w)=0$. 
\end{proposition}

A linear map $g$ on a vector space $V$ is called {\em
unipotent of index $n$} if $(g-I)^n(V)={\mathbf 0}$. 

\begin{lemma}~
\label{lem:phi}
\begin{enumerate}
\item $f_\ast:H_1(M,P;\Z) \to H_1(M,P;\Z)$ is unipotent of index $2$.
In particular, $\ker~\phi=\Fix(f_\ast^k)$ for all non-zero $k \in \Z$. 
\item $\phi\big(H_1(M,P;\Z)\big)$ is a submodule of $\langle \gamma_j \rangle$ of full rank.
Moreover, this rank is bounded from above by the genus of $M$. 
\item If $D\big(\Aff(M^\circ)\big)$ is non-elementary, then both $\hol
\circ \phi\big(H_1(M,P;\Z)\big)$ and $\hol(\ker~\phi)$ are
$\Z$-modules of 
rank $[k:\Q]$, where $k$ is the holonomy field. 
\end{enumerate}
\end{lemma}
\begin{proof}
We prove these statements in order. For all $x \in H_1(M,P;\Z)$, $\phi(x)$ is a linear combination of the 
$\{\gamma_j\}$. But, $i(\gamma_i,\gamma_j^\circ)=0$ for all $i$ and $j$. This implies statement (1). 

From equation (\ref{eq:phi}), we infer that
$\phi\big(H_1(M,P;\Z)\big) \subset \langle \gamma_j \rangle$.
Consider the map $\pi:H_1(M^\circ; \Z) \to H_1(M,P; \Z)$ induced by
the inclusion of $M^\circ \hookrightarrow M$. 
Define the map 
$$\eta:H_1(M,P; \Z) \to \langle\gamma^\circ_j \rangle, \ \ x \mapsto
\sum_{j} i(x, \gamma^\circ_j) t_j \gamma_j^\circ.$$  
Note that $\pi \circ \eta=\phi$.
We claim that the image of $\eta$ is a $\Z$-module of rank equal to
$\rk~ \langle\gamma^\circ_j\rangle$. 
If this is true, then the conclusion follows as
$\pi\left(\langle\gamma^\circ_j\rangle \right)=\langle\gamma_j\rangle$.
We now prove this claim. By non-degeneracy of 
$i:H_1(M, P; \Z) \times H_1(M^\circ; \Z) \to \Z$, it is equivalent to
show that if $x \in \ker(\eta)$ then $i(x, \gamma^\circ)=0$ for all
$\gamma^\circ \in \langle\gamma^\circ_j\rangle$.  
We will prove the contrapositive of this statement. Suppose $i(x,
\gamma^\circ) \neq 0$ for some $\gamma^\circ \in
\langle\gamma^\circ_j\rangle$. Then 
$i(x, \gamma^\circ_k) \neq 0$ for some $k$. We compute 
$$i\big(x, \eta(x)\big)=i\big(x,\sum_{j} i(x, \gamma^\circ_j) t_j
\gamma_j^\circ \big)=\sum_{j} t_j i(x, \gamma^\circ_j)^2.$$ 
Recall that each $t_j$ is non-zero and has the same sign. In addition,
$i(x, \gamma^\circ_k)\neq 0$, so $i\big(x, \phi(x)\big)\neq 0$. Therefore,
$\eta(x)\neq 0$. 

The inequality $\rk~ \langle\gamma_j \rangle \leq \mathit{genus}(M)$ follows from topology.
Note that the core curves of cylinders
are disjoint. Cutting along $g+1$ closed curves on a surface of genus
$g$ necessarily disconnects the surface. Hence, the maximal rank of
the span 
of the $\{\gamma_j\}$ is $\mathit{genus}(M)$, because the $\gamma_j$ have disjoint representatives.

Now we will consider statement (3). By applying an affine transformation to $M$, we may assume that
our cylinder decomposition is horizontal and that $\hol\big(H_1(M,P; \Z)\big)=k^2$. Then we have that 
$$D(f)=\left[\begin{array}{rr} 1 & \mu \\ 0 & 1 \end{array}\right],$$
with $0 \neq \mu \in k$.
Then for all $\gamma \in H_1(M,P;\Z)$, we have 
$$\hol \circ \phi(\gamma)=D(f) \hol(\gamma)-\hol(\gamma)=\big(\mu \hol_y(\gamma), 0\big),$$
where $\hol_y$ denotes the $y$-coordinate of the holonomy vector. We conclude
$$\rk~\hol \circ \phi\big(H_1(M,P;\Z)\big)=\rk~\hol_y\big(H_1(M,P;\Z)\big)=[k:\Q].$$
On the other hand, if $\gamma \in \ker~\phi$ then $\hol_y(\gamma)=0$. Thus
$\rk~\hol(\ker~\phi) \leq [k:\Q]$. Note that the fact that $\langle \gamma_j \rangle \subset \ker~\phi$ 
implies the opposite inequality, as $\rk~\hol\big( \langle \gamma_j \rangle \big)=[k:\Q]$.
This follows from 
Thurston's construction of pseudo-Anosov homeomorphisms, and assumes
that $D\big(\Aff(M^\circ)\big)$ is non-elementary. 
See Theorem 7 of \cite{T88}. 
\end{proof}

We first establish a corollary of statement (1) of the lemma.

\begin{corollary}
\label{cor:infinite_index}
Let $f \in \Aff(M^\circ)$ be a multi-twist, and let $w \in W$. If $f_\ast(w) \neq w$, then
$D\big(\Aff(\widetilde M_w, M^\circ)\big)$ is infinite index in $D\big(\Aff(M^\circ)\big)$.
\end{corollary}

Recall the definition of the holonomy-free subspace $W$ of
$H_1(M,P;\Z)$. Proposition \ref{prop:aff} stated that an 
element $f \in \Aff(M^\circ)$ lifted to an affine automorphism $\widetilde f \in \Aff(\widetilde M_w, M)$
if and only if $f_\ast(w)=\pm w$. 

\begin{theorem}[Lifting multi-twists]
\label{thm:multi-twist}
Assume $f \in \Aff(M^\circ)$ is a multi-twist and that $D\big(\Aff(M^\circ)\big)$ is non-elementary. 
Let the notation be as above. 
\begin{equation}
\label{eq: of the theorem}
\rk~W-\rk(W \cap \ker~\phi)=\rk~ \langle\gamma_j\rangle-[k:\Q] \leq
g-[k:\Q].
\end{equation}
In particular, $f_\ast$ acts trivially on $W$ if and only if 
$\rk~\langle\gamma_j\rangle=[k:\Q]$.
\end{theorem}

\begin{proof}
By linearity of $\phi$ and statement (2) of Lemma \ref{lem:phi}, 
$$\rk(\ker~\phi)=\rk~H_1(M,P;\Z)-\rk~\phi(H_1(M,P;\Z))=\rk~W+2[k:\Q]-\rk~\langle\gamma_j\rangle.$$
Now, note that $W \cap \ker~\phi=\ker~\hol |_{\ker~\phi}$. By linearity of $\hol$, we have
$$
\rk(\ker~\phi)=\rk (W \cap \ker~\phi)+\rk~\hol(\ker~\phi)=\rk (W \cap \ker~\phi)+[k:\Q],
$$
with the last equality following from statement (3) of the lemma.
Subtracting these two equations gives (\ref{eq: of the theorem}). 
The inequality follows from statement (2) of the lemma. 
\end{proof}

As an illustration of the use of Theorem \ref{thm:multi-twist}, we deduce:
\compat{Note strengthening of the corollary below.}

\begin{corollary}
\label{cor: improving HS}
Suppose $M$ is square-tiled and has cylinder decomposition in which all cylinders are homologous.
Then the Veech group of any recurrent $\Z$-cover is
Fuchsian of the first kind. 
\end{corollary}

\begin{proof}
In this case $\mathrm{rk} \langle \gamma_j \rangle =1$ and $k =\Q$, so
$f_* \in \ker~ \psi_0$. Since $Df_*$ is nontrivial, the result follows
from Theorem \ref{thm:kernel}. 
\end{proof}

\begin{remark}
In \cite[Theorem 2]{HS08}, Hubert and Schmith\"usen define a class of
$\Z$-covers of square tiled surfaces $\mathcal{O}^{\infty} \to
\mathcal{O}$. They show that if $\mathcal{O}$ has a one-cylinder
decomposition, then the Veech group of $\mathcal{O}^{\infty}$ is
Fuchsian of the first kind. Thus Corollary 
\ref{cor: improving HS} is an extension of the results of \cite{HS08}. 

\end{remark}
\section{Examples}

\subsection{Square tiled surfaces with homologous cylinders}
\compat{Notice the changes to this section.}
We give a construction which yields all square tiled surfaces with a horizontal cylinder decomposition consisting of
homologous cylinders.

Let $C_0, \ldots, C_{k-1}$ be cylinders all with the same rational circumference $c$, and each with rational width. For each $i=0, \ldots, k-1$ pick a rational interval exchange
of $T_i:[0,c) \to [0,c)$. Use $T_i$ to identify the bottom edge of $C_i$ to the top edge of $C_{i+1 (\textrm{mod } k)}$. Call the resulting surface $M$,
and let $P \subset M$ be a finite set of points with rational coordinates. Then there is a horizontal cylinder decomposition of $M^\circ$, all of whose cylinders
are homologous. So, by Corollary \ref{cor: improving HS}, any recurrent $\Z$-cover of has a Veech group which is Fuchsian of the first kind. 

\def\eW{\mathit{W}}

\begin{figure}
\includegraphics{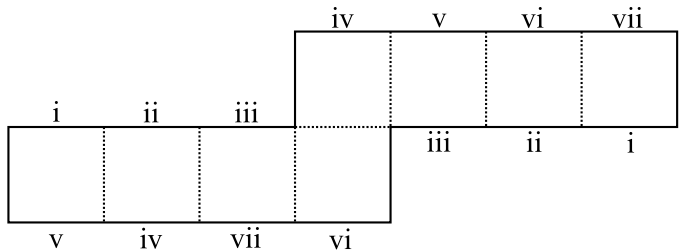}
\caption{The eierlegende Wollmilchsau surface. Horizontal edges are
glued together as indicated by the roman numerals. Vertical edges are
glued to their opposite (by horizontal translations).} 
\label{fig:eW}
\end{figure}

The term {\em eierlegende Wollmilchsau} refers to the square tiled
surface, $\eW$, whose properties were first studied by Herrlich and
Schmith\"usen 
\cite{HS08}. It can be obtained by the above construction. 
See figure \ref{fig:eW}. 
This is a surface of genus three with four cone singularities, each with cone angle $4 \pi$. 
Let $P$ denote the set of these singularities.
The Veech group of $\eW^\circ$ is $\widehat
\SL(2,\Z)$, the group of integer matrices of determinant $\pm 1$. 

\begin{proposition}
Any $\Z$-cover of $\eW^\circ$ has a Veech group that contains the congruence $4$ subgroup of
$\SL(2, \Z)$.
\end{proposition}
\begin{proof}
The horizontal direction has a multi-twist $\phi$ in a pair of homologous cylinders with derivative 
$D(\phi)=\left[\begin{array}{rr} 1 & 4 \\ 0 & 1\end{array}\right]$. For any $B \in \widehat\SL(2,\Z)=\Gamma(\eW^\circ)$,
there is a multi-twist $\phi_B$ in a pair of homologous cylinders with derivative $D(\phi_B)=BD(\phi)B^{-1}$. By Corollary \ref{cor: improving HS},
each $\phi_B$ lifts to any recurrent $\Z$-cover. The derivatives of these elements generate the congruence $4$ subgroup of
$\SL(2, \Z)$.
\end{proof}

\subsection{A question of Hubert and Schmith\"usen}
\label{sect:HS}
\begin{figure}
\includegraphics[height=1.5in]{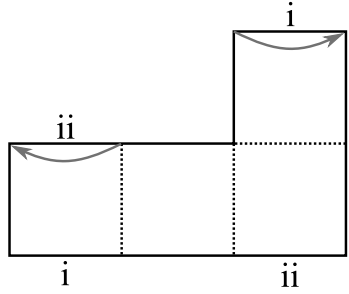}
\caption{The surface $Z_{3,1}$ and the cycle $w$.}
\label{fig:z31}
\end{figure}

We consider a surface defined in
\cite{HS08}. Let $Z_{3,1}$ be as in figure \ref{fig:z31}, let $w$ be the cycle marked on figure \ref{fig:z31} and let
$Z_{3,1}^{\infty}$ be the corresponding $\Z$-cover. Since $\hol(w)=0$
this is a recurrent $\Z$-cover. Hubert and Schmith\"usen proved that
the Veech group of $Z_{3,1}$ is not a lattice, but, since the genus of
$Z_{3,1}$ is 2, $\rk~W_0 =2$ so Theorem \ref{thm:dimension_2} implies that the Veech
group of $Z_{3,1}^{\infty}$ is Fuchsian of the first kind. This
answers a question raised in \cite{HS08}. 

\subsection{A double cover of the octagon}
\def\oct{{O}}
Let $X$ denote the polygon shown on the left side of figure \ref{fig:oct}. 
The translation surface $\oct$ is obtained by applying the Zemlyakov-Katok
unfolding construction to $X$ \cite{ZK}. The surface $\oct$ is a double cover of the regular octagon 
with opposite sides identified, as depicted on the right side of figure \ref{fig:oct}. 
The surface $\oct$ is of genus
$3$ with two cone singularities, each with cone angle $6 \pi$. 

\combarak{Pascal told me that Veech had already considered a double
cover of the octagon as well, see \S4 of his paper. I am guessing it
is this cover although it would be better if you look, I have a hard
time with Veech's terminology}

\compat{The double cover considered by Veech is different. I added a description in terms of a billiard table which may make this clearer.}

\begin{figure}
\includegraphics[width=2.5in]{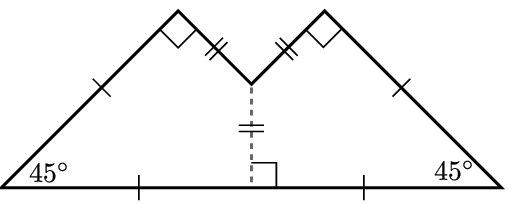}
\hspace{0.5in}
\includegraphics[width=3in]{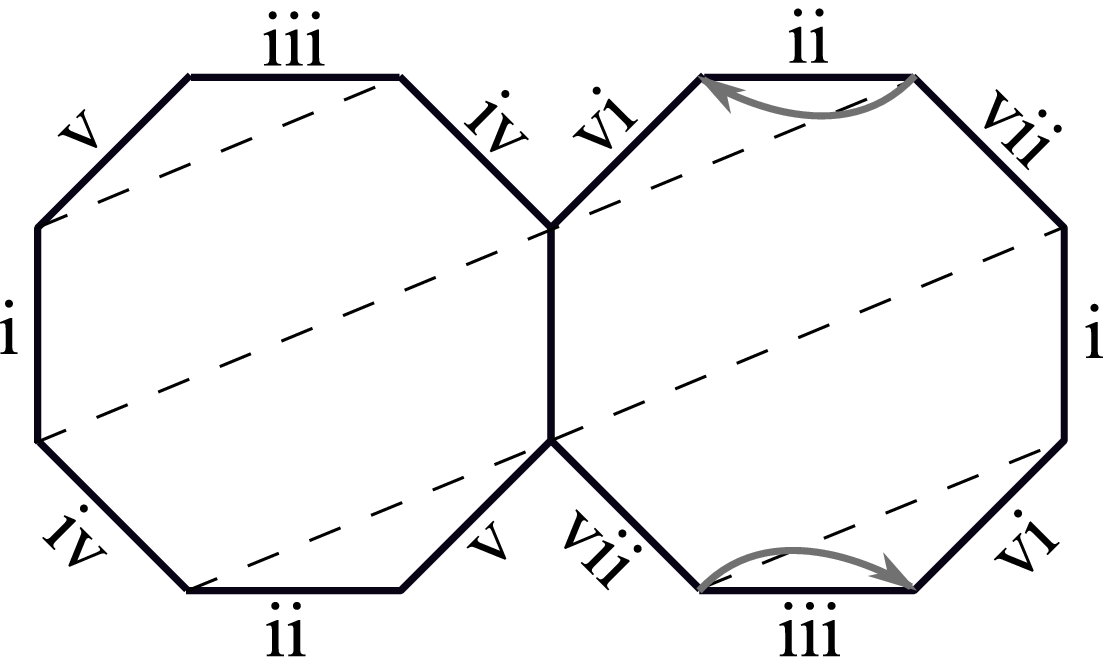}
\caption{The polygon $X$ and the surface $\oct$.}
\label{fig:oct}
\end{figure}

Let $P$ consist of the two singularities of $\oct$.
The orientation preserving part of the Veech group is generated by the
derivatives of the following affine automorphisms. 
\begin{itemize}
\item $h \in \Aff(\oct^\circ)$ is the right multi-twist in the
horizontal cylinder decomposition. We have
$D(h)=\left[\begin{array}{rr} 1 & 2 +\sqrt{2} \\ 0 & 1
\end{array}\right].$ 
\item $g \in \Aff(\oct^\circ)$ is the right multi-twist in the
cylinder decomposition in the direction of angle $\pi/4$. We have  
$D(g)=\left[\begin{array}{rr} -\sqrt{2} & 1 + \sqrt{2}\\ 
-1-\sqrt{2} &  2 + \sqrt{2} \end{array}\right].$
\item $f \in \Aff(\oct^\circ)$ is the right multi-twist in the
cylinder decomposition in the direction of angle $\pi/8$.  
$D(f)=\left[\begin{array}{rr} -1-\sqrt{2} & 4 + 3 \sqrt{2}\\ 
-\sqrt{2} &  3 + \sqrt{2} \end{array}\right].$
\item The two elements in $\Aff(\oct^\circ)$ with derivative $-I$. 
\end{itemize}
The orientation preserving part of the Veech group
$D\big(\Aff(\oct^\circ))$ is an index two subgroup of a $(4, \infty,
\infty)$-triangle group.

\begin{proposition}\label{prop: octagon}
For any $w \in W \subset H_1(\oct, P; \Z)$, there is a lift of $f \in \Aff(\oct^\circ)$ to 
$D\big(\Aff(\widetilde \oct_w, \oct^\circ)\big)$. In particular, 
$D\big(\Aff(\widetilde \oct_w, \oct^\circ)\big)$ is always a Fuchsian group of the first kind.
\end{proposition}
\begin{proof}
The affine automorphism $f$ is a multi-twist which preserves a
cylinder decomposition consisting of two cylinders.  
By the multi-twist theorem, it fixes all of $W$. By Theorem
\ref{thm:kernel}, $D\big(\Aff(\widetilde \oct_w)\big)$ is a Fuchsian group
of the first kind. 
\end{proof}

The following gives an example of an infinite translation surface with
non-arithmetic Veech group. 

\begin{proposition}
There exists a $w_1 \in W$ for which $D\big(\Aff(\widetilde \oct_{w_1},
\oct^\circ)\big)$ is an infinitely generated 
Fuchsian group of the first kind, and a $w_2 \in W$ for which 
$D\big(\Aff(\widetilde \oct_{w_2}, \oct^\circ)\big)$ contains the 
lattice $\langle D(f), D(g), D(h) \rangle \subset D\big(\Aff(\oct^\circ)\big)$.
\end{proposition}
\begin{proof}
We saw in the previous proposition that $f$ always lifts. As $\oct$ is
genus $3$, the multi-twist theorem implies that 
$\Fix_W(g_\ast)$ and $\Fix_W(h_\ast)$ are at worst codimension $1$ inside $W$. Note that
$\dim~W=3$. Thus, we can find a non-zero $w_2 \in \Fix_W(g_\ast) \cap \Fix_W(h_\ast)$. 
As
$D\big(\Aff(\oct^\circ)\big)$ is generated by $\langle D(f), D(g), D(h)\rangle$, we see
$D\big(\Aff(\widetilde \oct_{w_2}, \oct^\circ)\big)=D\big(\Aff(\oct^\circ)\big)$.

To see that there is a $w_1 \in W$ for which $D\big(\Aff(\widetilde
\oct_{w_1}, \oct^\circ)\big)$ is infinitely generated, 
it is sufficient to show that $D\big(\Aff(\widetilde \oct_{w_2},
\oct^\circ)\big)$ is infinite index in $D\big(\Aff(\oct^\circ)\big)$. 
By Corollary \ref{cor:infinite_index} and the multi-twist theorem, it
is sufficient to check that the span of the core curves of a cylinder
decomposition  
span a rank three submodule of $H_1(\oct,P;\Z)$. This is true for both the
horizontal direction and the direction of angle $\pi/4$. 
\end{proof}

It turns out that there is only one non-zero $w \in W$ up to scaling which is fixed by
$f_\ast$, $g_\ast$ and $h_\ast$. This $w$ is the homology class shown
in grey in 
figure \ref{fig:oct}.

\vspace{1em}
\noindent
{\bf Acknowledgements.} The first author would like to thank the
Institute for Advanced studies in Mathematics at Ben-Gurion University
of the Negev 
for funding a trip to Israel during which most of the ideas in this
paper were worked out. The authors would like to thank Yair Minsky and
Pete Storm for ideas which went into the proof 
of Theorem \ref{thm:dimension_2}.

\bibliographystyle{amsalpha}
\bibliography{bibliography}
\end{document}